\documentclass{amsart}

\usepackage{amsfonts}
\usepackage{amssymb}
\usepackage{amsmath}
\usepackage{graphicx}
\usepackage{amsthm}
\usepackage{xcolor}
\usepackage{hyperref}
\usepackage{nicematrix}

\def\({\left(}
\def\){\right)}

\numberwithin{equation}{section}

\newtheorem{theorem}{Theorem}[section]
\newtheorem*{theorem*}{Remark}
\newtheorem{definition}[theorem]{Definition}
\newtheorem{proposition}[theorem]{Proposition}
\newtheorem{lemma}[theorem]{Lemma}
\newtheorem{remark}[theorem]{Remark}

\newtheorem{corollary}[theorem]{Corollary}
\newtheorem{conjecture}[theorem]{Conjecture}

\makeatletter
\@namedef{subjclassname@2010}{%
  \textup{2010} Mathematics Subject Classification}
\makeatother

\begin{document}

\title{Expanders on matrices over a finite chain ring, ${\rm I}$}

\author{Dung M. Ha}
\address{School of Applied Mathematics and Informatics, Hanoi University of Science and Technology,
Hanoi, Vietnam}
\email{dung.hm200096@sis.hust.edu.vn}

\author{Hieu T. Ngo}
\address{Institute of Mathematics, 
Vietnam Academy of Science and Technology, Hanoi, Vietnam}
\email{nthieu@math.ac.vn}

\thanks{}

\begin{abstract} 
In this work and its sequel, we study the expanding phenomenon of matrices over a finite chain ring of large residue field.
A sum-product estimate is proved. 
It is showed that $x+yz$ is a moderate expander on $n\times n$ matrices with exponent $\frac{n+1}{6}$.
These results generalise the main theorems in the recent work \cite{XG} of Xie and Ge. The proofs use spectral graph theory and elementary divisor theory.
\end{abstract}

\maketitle

\section{Introduction}


In 1983, Paul Erd\"os and Endre Szemer\'edi \cite{ES83} proposed the Sum-Product Conjecture, which still remains open today:
\begin{conjecture}(Erd\"os-Szemer\'edi)
For every $\epsilon>0$, there is $n_0=n_0(\epsilon)$ such that for any $A\subset \mathbb{N}$ with $|A|>n_0$, one has
$$
|(A+A) \cup (A\cdot A)| > |A|^{2-\epsilon}. 
$$
\end{conjecture}
\noindent In plain words, this sum-product phenomenon indicates that the sum set and the product set of a set of integers cannot both be small compared to $|A|$. This problem has a long and exciting development. In place of the exponent $2-\epsilon$, Erd\"os and Szemer\'edi \cite{ES83} established this inequality with the exponent $1+\delta$ for some $\delta>0$, and the best current record is the exponent $\frac43 + \frac{2}{1167}$ by Rudnev and Stevens \cite{RS21}.

In 2004, Bourgain, Katz and Tao \cite{BKT} initiated the study of the finite field analogue of the Sum-Product Problem and established the following theorem: for a subset $A$ of the prime field $\mathbb{F}_p$, if $p^\delta<|A|<p^{1-\delta}$ for some $\delta>0$, then $\max(|A+A|,|A\cdot A|)\gg |A|^{1+\kappa}$ for some $\kappa=\kappa(\delta)>0$. Here the asymptotic notation $X \gg Y$ means that there exists some constant $c>0$ such that $X \geq cY$;
if $X=X(n)$ and $Y=Y(n)$ depend on some parameter $n$, the asymptotic $X \gg_n Y$ means that there exists some constant $c(n)>0$ such that $X \geq c(n)Y$.
Since the work \cite{BKT}, there have been many interesting findings on the sum-product phenomenon over a finite field or a modular ring. Around 2007, V.~Vu \cite{Vu08} and V.~Le \cite{Vinh11} independently introduced the application of spectral graph theory to the Sum-Product Problem in a finite setting, besides tools from incidence geometry, Fourier analysis, and exponential sums. The new toolkit involving spectra of graphs provides not only new proofs but also improvements of known estimates as well as many new results (see, e.g., \cite{Solymosi09}).

The sum-product phenomenon can be bootstrapped to a multivariate question: given a finite subset $A$ of a ring $S$, under which multivariable function $f:S^d\to S$ does the image set $f(A,\dots,A)$ grow large? 
In particular, detecting this expansive phenomenon in finite rings, commutative or non-commutative, is a topic that has attracted wide interest in the last two decades. 

Let $\mathbb{F}_{q}$ denote the finite field of order $q$.
Over finite fields, Hart, Li, Shen \cite{HLS13} and Tao \cite{Tao15} defined several classes of expanding functions. 

\begin{definition}[Expanders on finite fields]
Let $f: \mathbb{F}_{q}^d\to \mathbb{F}_{q}$ be a function.
\begin{enumerate}
    \item We say that $f$ is a \emph{moderate expander} on $\mathbb{F}_{q}$ with exponent $\kappa>0$ if there is an absolute constant $C>0$ such that for any subset $A\subset \mathbb{F}_{q}$ with $|A| \geq Cq^{1-\kappa}$, one has $|f(A,\dots,A)| \geq C^{-1}q$.
    \item We say that $f$ is a \emph{strong expander} on $\mathbb{F}_{q}$ with exponent $\kappa>0$ if there are absolute constants $C,C'>0$ such that for any subset $A\subset \mathbb{F}_{q}$ with $|A| \geq Cq^{1-\kappa}$, one has $|f(A,\dots,A)| \geq q-C'$.
    \item We say that $f$ is a \emph{very strong expander} on $\mathbb{F}_{q}$ with exponent $\kappa>0$ if there is an absolute constant $C>0$ such that for any subset $A\subset \mathbb{F}_{q}$ with $|A| \geq Cq^{1-\kappa}$, one has $f(A,\dots,A)=\mathbb{F}_{q}$.
\end{enumerate}
\end{definition}

There have been found many interesting polynomial expanders with exponents on finite fields. 

\begin{theorem} On $\mathbb{F}_q$, one has the following polynomial expanders.
\begin{enumerate}
    \item (Hart, Iosevich, Solymosi \cite[Theorem 1.4]{HIS07}) $(x-y)(z-t)$ is a very strong expander with exponent $\kappa=\frac14$.
    \item (Bennett, Hart, Iosevich, Pakianathan, Rudnev \cite[Corollary 1.8]{BHIPR17}) $(x-y)(z-t)$ is a moderate expander with exponent $\frac13$. 
    \item (Hart, Iosevich \cite[Theorem 1.1]{HI08}) $xy+zt$ is a strong expander with exponent $\frac14$ and $C=1$, and is a moderate expander with exponent $\frac13$.
    \item (Shparlinski \cite[Equation $(9)$]{Shparlinski08}) $x+yz$ is a moderate expander with exponent $\frac13$.
    \item (V.~Le \cite{Vinh13}) $x(y+z)$ and $x+(y-z)^2$ is a moderate expander with exponent $\frac13$.
    \item (V.~Le \cite{Vinh13}) $(x-y)^2+zt$ is a moderate expander with exponent $\frac38$.
\end{enumerate}
\end{theorem}

A commutative analogue of a finite field is a modular ring of prime power order. In commutative ring theory, there is a notion that captures both these rings as special cases:
a \emph{finite chain ring} is a finite commutative ring which is also a local principal ideal ring (see Section \ref{subsection:fcr}).
Let $R$ be a finite chain ring. 
The ring $R$ has a unique maximal ideal $\mathfrak{m}$, which is a principal ideal. The residue field $k=R/\mathfrak{m}$ is a finite field $\mathbb{F}_q$, say. Let $r$ denote the nilpotency of $\mathfrak{m}$, i.e.~$r$ is the smallest positive integer $t$ such that $\mathfrak{m}^t=0$.
Interesting examples of a finite chain ring include:
\begin{itemize}
    \item finite fields;
    \item $\mathbb{Z}/p^r\mathbb{Z}$ where $p$ is a prime number;
    \item $\mathcal{O}/(\pi^r)$ where $\mathcal{O}$ is the ring of integers in a number field and $\pi$ is a prime in $\mathcal{O}$;
    \item $\mathbb{F}_q[x]/(f^r)$ where $f$ is an irreducible polynomial in $\mathbb{F}_q[x]$.
\end{itemize}
In recent literature, some authors use the term ``finite valuation ring'' with the same meaning as a finite chain ring. 
However, we choose the term ``finite chain ring'' instead 
(see Section \ref{subsection:fcr}).

We will consider the expanding phenomenon on finite chain rings of large residue field:

\begin{definition}[Expanders on finite chain rings]\label{definition:expander-fcr}
Let $(R,\mathfrak{m})$ be a finite chain ring. Suppose that the residue field $k=R/\mathfrak{m}$ has cardinality $q$ and that $R$ has cardinality $q^r$. Let $f: R^d\to R$ be a function.
\begin{enumerate}
    \item We say that $f$ is a \emph{moderate expander} on $R$ with exponent $\kappa>0$ if if there is an absolute constant $C>0$ such that for any subset $A\subset R$ with $|A| \geq Cq^{r-\kappa}$, one has $|f(A,\dots,A)| \gg q^r$.
    \item We say that $f$ is a \emph{strong expander} on $R$ with exponent $\kappa>0$ if if there are absolute constant $C,C'>0$ such that for any subset $A\subset R$ with $|A| \geq Cq^{r-\kappa}$, one has $|f(A,\dots,A)| \geq q^r-C'$.
    \item We say that $f$ is a \emph{very strong expander} on $R$ with exponent $\kappa>0$ if if there is an absolute constant $C>0$ such that for any subset $A\subset R$ with $|A| \geq Cq^{r-\kappa}$, one has $f(A,\dots,A)=R$.
\end{enumerate}
\end{definition}

\begin{remark}
On taking $r=1$, we recover expanders on finite fields.
\end{remark}

Recent researches reveal an abundant supply of polynomial expanders with exponents on finite chain rings. 

\begin{theorem}
Let $(R,\mathfrak{m})$ be a finite chain ring. Suppose that the residue field $k=R/\mathfrak{m}$ has cardinality $q$ and that $R$ has cardinality $q^r$. 
On $R$, one has the following polynomial expanders.
\begin{enumerate}
    \item (Yazici \cite[Theorem 1.1]{Yazici}) $x+yz$ is a moderate expander with exponent $\frac{1}{3}$.
    \item (Dao, Koh, Le, Mirzaei, Mojarrad, Pham \cite[Theorem 1.5]{DKea}) $x(x+y)z+t$, $x(x+y)+zt$, $x(x+y)(z+t)$, $z(x(x+y)+t)$, $(x(x+y)-z)^{2}+t$ and $(z-t)^{2}+x(x+y)$ are moderate expanders with exponent $\frac{3}{8}$.
\end{enumerate}
\end{theorem}

A theme of considerable interest is to study sums and products in matrix spaces, starting with Chang \cite{Chang07}. As noted in \cite{Chang07}, one interesting point about this non-commutative quest is that there are large sets of matrices whose sum sets and product sets do not grow.
Many additive combinatorics problems have been generalised from commutative rings to matrices over commutative rings, albeit with different techniques.
Therefore, it is curious to study the expanding phenomenon on matrices.
For a ring $S$ with $1$, let ${\rm M}_n(S)$ (resp.~${\rm GL}_n(S)$) denote the ring of $n\times n$ matrices (resp.~the group of invertible $n\times n$ matrices) with entries in $S$.

\begin{definition}[Expanders on matrix algebras over finite fields]
Let $f:{\rm M}_n(\mathbb{F}_{q})^d\to {\rm M}_n(\mathbb{F}_{q})$ be a function.
\begin{enumerate}
    \item We say that $f$ is a \emph{moderate expander} on ${\rm M}_n(\mathbb{F}_{q})$ with exponent $\kappa(n)>0$ if there is a constant $C(n)>0$ such that for any subset $A\subset {\rm M}_n(\mathbb{F}_{q})$ with $|A| \geq C(n) q^{n^{2}-\kappa(n)}$, one has $|f(A,\dots,A)| \geq C(n)^{-1} q^{n^2}$.
    \item We say that $f$ is a \emph{strong expander} on ${\rm M}_n(\mathbb{F}_{q})$ with exponent $\kappa(n)>0$ if there are constants $C(n),C'(n)>0$ such that for any subset $A\subset {\rm M}_n(\mathbb{F}_{q})$ with $|A| \geq C(n) q^{n^{2}-\kappa(n)}$, one has $|f(A,\dots,A)| \geq q^{n^2}-C'(n)$.
    \item We say that $f$ is a \emph{very strong expander} on ${\rm M}_n(\mathbb{F}_{q})$ with exponent $\kappa(n)>0$ if there is a constant $C(n)>0$ such that for any subset $A\subset {\rm M}_n(\mathbb{F}_{q})$ with $|A| \geq C(n) q^{n^{2}-\kappa(n)}$, one has $f(A,\dots,A)={\rm M}_n(\mathbb{F}_{q})$.
\end{enumerate}
\end{definition}

Several noncommutative polynomial expanders over matrix algebras over a finite field have been produced, although not as many as the cases of finite fields and finite chain rings.

\begin{theorem}\label{theorem:expanders-matrix-finite-field} 
On ${\rm M}_n(\mathbb{F}_q)$, one has the following noncommutative polynomial expanders.
\begin{enumerate}
    \item (V.~Le, T.~Nguyen \cite[Theorem 1.11]{TheVinh20}) $xy+z+t$ is a very strong expander with exponent $\kappa=\frac{1}{4}$.
    \item (V.~Le, T.~Nguyen \cite[Corollary 1.8]{TheVinh20}) $x(y+z)$ is a moderate expander with exponent $1$.
    \item (V.~Le, T.~Nguyen \cite[Theorem 1.10]{TheVinh20}) $x+yz$ is a moderate expander with exponent $\frac{1}{3}$.
    \item (Xie, Ge \cite{XG}) $x+yz$ is a moderate expander with exponent $\frac{n+1}{6}$.
\end{enumerate}
\end{theorem}

The moderate expander $x+yz$ above is a consequence of a more general estimate:

\begin{theorem}\label{theorem-XieGe:expander}\cite[Theorem 1.5]{XG}
For every positive integer $n$, there exists a constant $C(n)>0$ such that 
if $A,B,C\subset {\rm M}_n(\mathbb{F}_q)$, then
$$
|A+BC| \geq C(n) \min \left( q^{n^{2}}, \frac{|A||B||C|}{q^{2n^{2}-\frac{n+1}{2}}} \right).
$$
\end{theorem} 

In the same work \cite{XG}, the authors also established a new sum-product estimate on matrix algebras over finite fields.

\begin{theorem}\label{theorem-XieGe:sumproduct}\cite[Theorem 1.6]{XG}
For every positive integer $n$, there exists a constant $C(n)>0$ such that for each $A\subset {\rm M}_n(\mathbb{F}_q)$ satisfying $|A|\geq C(n)q^{n^{2}-1}$, 
$$
\max ( |A+A|,|AA| )
    \geq C(n)^{-1} \min \left( \frac{|A|^{2}}{q^{n^{2}-\frac{n+1}{4}}} , q^{\frac{n^{2}}{3}}|A|^{\frac{2}{3}} \right) .
$$
\end{theorem}

It is interesting to ask whether the sum-product phenomenon and the expanding phenomenon happen in matrix algebras over finite chain rings. It seems that this question has not been studied in the literature. We first define expander classes in this context.

\begin{definition}[Expanders on matrix algebras over finite chain rings]
Let $(R,\mathfrak{m})$ be a finite chain ring. Suppose that the residue field $k=R/\mathfrak{m}$ has cardinality $q$ and that $R$ has cardinality $q^r$. 
Let $f:{\rm M}_n(R)^d\to {\rm M}_n(R)$ be a function.
\begin{enumerate}
    \item We say that $f$ is a \emph{moderate expander} on ${\rm M}_n(R)$ with exponent $\kappa(n)>0$ if there is a constant $C(n)>0$ such that for any subset $A\subset {\rm M}_n(R)$ with $|A| \geq C(n) q^{rn^{2}-\kappa(n)}$, one has $|f(A,\dots,A)| \geq C(n)^{-1} q^{rn^2}$.
    \item We say that $f$ is a \emph{strong expander} on ${\rm M}_n(R)$ with exponent $\kappa(n)>0$ if there are constants $C(n),C'(n)>0$ such that for any subset $A\subset {\rm M}_n(R)$ with $|A| \geq C(n) q^{rn^{2}-\kappa(n)}$, one has $|f(A,\dots,A)| \geq q^{rn^2}-C'(n)$.
    \item We say that $f$ is a \emph{very strong expander} on ${\rm M}_n(R)$ with exponent $\kappa(n)>0$ if there is a constant $C(n)>0$ such that for any subset $A\subset {\rm M}_n(R)$ with $|A| \geq C(n) q^{rn^{2}-\kappa(n)}$, one has $f(A,\dots,A)={\rm M}_n(R)$.
\end{enumerate}
\end{definition}

In this paper, we generalize the two main theorems of \cite{XG}.

\begin{theorem}\label{theorem-result:expander}
Let $(R,\mathfrak{m})$ be a finite chain ring. Suppose that the residue field $k=R/\mathfrak{m}$ has cardinality $q$ and $R$ has cardinality $q^r$.
If $A,B,C\subset {\rm M}_n(R)$, then for every positive integer $n$ there exists a constant $C(n)>0$ such that
$$
|A+BC| \geq C(n)\min \left( q^{rn^{2}}, \frac{|A||B||C|}{q^{2rn^2-\frac{n+1}{2}}} \right).
$$
\end{theorem} 

\begin{theorem}\label{theorem-result:sumproduct}
Let $(R,\mathfrak{m})$ be a finite chain ring. Suppose that the residue field $k=R/\mathfrak{m}$ has cardinality $q$ and $R$ has cardinality $q^r$.
For every positive integer $n$, there exists a constant $C(n)>0$ such that the following holds. If $A\subset {\rm M}_n(R)$ satisfies $|A|\geq C(n)q^{rn^{2}-1}$, then
$$
\max ( |A+A|,|AA| )
    \geq C(n)^{-1} \min \left( \frac{|A|^{2}}{q^{rn^{2}-\frac{n+1}{4}}} , q^{\frac{rn^{2}}{3}}|A|^{\frac{2}{3}} \right) .
$$
\end{theorem}

\begin{remark}
    In Theorems \ref{theorem-result:expander} and \ref{theorem-result:sumproduct}, on specialising $R=\mathbb{F}_q$, whence $r=1$, one recovers Theorems \ref{theorem-XieGe:expander} and \ref{theorem-XieGe:sumproduct} respectively.
\end{remark}

We note an immediate consequence of Theorem \ref{theorem-result:expander}.
\begin{corollary}\label{corollary-result:expander}
If $A,B,C\subset {\rm M}_n(R)$ satisfy $|A||B||C| \gg_n q^{3rn^2-\frac{n+1}{2}}$, then $|A+BC| \gg_n q^{rn^{2}}$.
In particular, $x+yz$ is a moderate expander on ${\rm M}_n(R)$ with exponent $\frac{n+1}{6}$.
\end{corollary} 

Our results offer first examples of a sum-product estimate and a noncommutative polynomial expander on matrix algebras over finite chain rings.
The strategy of proofs follows a standard approach: we first construct graphs from the original combinatorial problems, and then apply theorems concerning spectra of graphs to count solutions to matrix equations. Many authors have adopted this line of arguments successfully \cite{Vu08,Vinh11,ThangVinh15,ThangVinh17,XG,MPW21}.
The generalisation from matrix algebras over finite fields to matrix algebras over finite chain rings, is, however, not obvious. For this purpose, we first observe that a finite chain ring is an elementary divisor ring, whence matrices over a finite chain ring possess Smith normal forms. Then we apply Smith normal form as a crucial ingredient to construct graphs and to estimate the number of solutions of a certain matrix equation (Proposition \ref{proposition:count-matrix-solutions}). On combining spectral graph theory with linear algebra over finite fields and finite chain rings, we are able to generalize the new results in \cite{XG}. 

The paper is organised as follows.
In Section \ref{section:preliminaries}, we recall some facts concerning finite chain rings, eigenvalues of graphs, and Smith normal form.
In Section \ref{section:counting-fcr}, we count matrices with a prescribed Smith normal form and count solutions of a system of linear equations.
In Section \ref{section:matrix-equation}, we study the matrix equation $ab+ef=c+d$ and prove the key estimate (Proposition \ref{proposition:count-matrix-solutions}) on the number of its solutions. 
Finally, Section \ref{section:proofs-of main-theorems} contains the proofs of the main results (Theorems \ref{theorem-result:expander} and \ref{theorem-result:sumproduct}).

\section{Preliminaries}\label{section:preliminaries}

\subsection{Finite chain rings}\label{subsection:fcr}

In this paragraph, we review the notion of a finite chain ring.

\begin{definition} Let $S$ be a commutative, unital ring.
\begin{enumerate}
    \item $S$ is called a \emph{local ring} if it has only one maximal ideal.
    \item $S$ is called a \emph{principal ideal ring} if all of its ideals are principal. 
    \item $S$ is called a \emph{chain ring} if its ideals form a chain under set-theoretic inclusion. 
\end{enumerate}
\end{definition}

\begin{proposition}\label{proposition:fcr-characterization}\cite[Proposition 2.1]{DLP04}
Let $S$ be a finite commutative ring.
The following are equivalent:
\begin{enumerate}
    \item $S$ is a chain ring.
    \item $S$ is a local ring and its maximal ideal is principal.
    \item $S$ is a local principal ideal ring.
\end{enumerate}
\end{proposition}

A finite commutative ring satisfying the equivalent conditions of Proposition \ref{proposition:fcr-characterization} is called a \emph{finite chain ring}. 
In Coding Theory literature, finite chain rings find many applications \cite{DLPS09}. 

Throughout this paper, we shall let $R$ denote a finite chain ring. Its unique maximal ideal is denoted by $\mathfrak{m}$. We let $\pi$ be a uniformizer, so $\mathfrak{m}=(\pi)$. The residue field $k=R/\mathfrak{m}$ is a finite field of order $q$. We write $r$ for the nilpotency of $\pi$, i.e.~$r$ is the smallest positive integer $t$ such that $\pi^t=0$. 
We write ${\rm val}$ for the natural valuation over $R$; this is defined by ${\rm val}(0)=r$ and for $x\in R\backslash\{0\}$, ${\rm val}(x)=k$ if and only if $x=u\pi^{k}$ for some unit $u\in R^\times$.

The following fact is standard (see for instance \cite{ThangVinh17,HamThangVinh17,Hiep21}).

\begin{lemma}
\label{lemma:prelim-fcr-size}
Suppose that $k\in \mathbb{N}$ satisfies $k\leq r$ and that $a\in R$. Then
$$
|a+\mathfrak{m}^{k}|=|\mathfrak{m}^{k}| = q^{r-k} .
$$
\end{lemma}

\subsection{Spectral graph theory}

This paragraph recalls some well-known results on eigenvalues of (undirected) graphs.
We say that a graph is \emph{regular of degree $d$} if each of its vertices is adjacent to $d$ other vertices.
We say that a bipartite graph is \emph{biregular} of degree $d$ if, considered as a usual graph, it is regular of degree $d$.

\begin{lemma}\label{lemma:prelim-regular-graph-largest-eval}
Let $G$ be a regular graph of degree $d$ and let $A_{G}$ be the adjacency matrix of $G$. Then $d$ is the eigenvalue of $A_{G}$ with the largest absolute value. In particular, if $\lambda$ is an eigenvalue of $A_G$, then $|\lambda| \leq d$.
\end{lemma}

Let $G$ be a bipartite graph and $A_G$ be its adjacency matrix. 
We order the eigenvalues of $A_G$ such that 
$|\lambda_1| \geq |\lambda_2| \geq \dots \geq |\lambda_n|$.
Since $G$ is bipartite, the eigenvalues of $A_G$ are symmetric around the origin, so we have $\lambda_1+\lambda_2=0$. 
We recall two well-known facts concerning a biregular bipartite graph.
Lemma \ref{lemma:expander-mixing} is known as the Expander Mixing Lemma; its proof can be found, for instance, in \cite[Lemma 2.1]{LNP19}.

\begin{lemma}[Expander Mixing Lemma]\label{lemma:expander-mixing}
Let $G$ be a biregular bipartite graph with parts $U$ and $V$. Then for every pair $X\subset U$ and $Y\subset V$, the number of edges between $X$ and $Y$, denoted by $e(X,Y)$, satisfies
$$
\left| e(X,Y)-\frac{\deg (U)}{|V|} \right| \leq |\lambda_{3}|\sqrt{|X||Y|}
$$
where $\lambda_{3}$ is the eigenvalue with the third largest absolute value of $A_{G}$.
\end{lemma}

\begin{lemma}[\cite{Vinh11}]\label{lemma:prelim-biregular-eval}
Let $G$ be a biregular bipartite graph whose parts $U$ and $V$ have sizes  $|U|=m$ and $|V|=n$. We label vertices of $G$ from $1$ to $m+n$. Then 
$$
A_{G}=
\begin{pNiceMatrix}
0 & N \\
N^{T} & 0
\end{pNiceMatrix}
$$
where $N$ is the $m \times n$ matrix and $N_{ij}=1$ if and only if there is an edge between $i$ and $j$. Let $w^{(3)}=(u_{1},\dots,u_{m},v_{1},\dots,v_{n})^{T}$ be an eigenvector of $A_{G}$ corresponding to the eigenvalue $\lambda_{3}$. Then:
\begin{enumerate}
\item $(u_{1},\dots,u_{m})^{T}$ is an eigenvector of $NN^{T}$ corresponding to the eigenvalue $\lambda_{3}^{2}$, and
\item $J(u_{1},\dots,u_{m})^{T}=0$ where $J$ is the $m\times m$ all-ones matrix.
\end{enumerate}
\end{lemma}

We will also need a classical fact about the largest eigenvalue of a real symmetric matrix (cf. \cite[Section 1]{Fulton}, \cite{KT}).
\begin{lemma}
\label{lemma:prelim-sum-symmetric-matrices-eigenvalue}
For a real symmetric matrix $M$, let $\lambda_{max}(M)$ denote the eigenvalue of $M$ with the largest absolute value. If $A$ and $B$ are two real symmetric matrices, then
$$|\lambda_{max}(A+B)| \leq 
    |\lambda_{max}(A)| + |\lambda_{max}(B)|.$$
\end{lemma}

\subsection{Elementary divisor theory}

In this paragraph, we review the theory of elementary divisors and normalise the Smith normal form for our computations.

\begin{definition}\cite{Kaplansky},\cite[Definitions 15.7 and 15.10]{Brown}
Let $S$ be a commutative, unital ring.
Let ${\rm diag}(d_1,d_2,\dots)$ denote a (possibly rectangular) matrix with $d_1,d_2,\dots$ on the diagonal and with zeros on the other entries. 
A matrix $A\in {\rm M}_{m\times n}(S)$ is said to \emph{admit diagonal reduction} if there exist invertible matrices $P\in {\rm GL}_m(S)$ and $Q\in {\rm GL}_n(S)$ such that 
$$
PAQ={\rm diag}(d_1,d_2,\dots) \in {\rm M}_{m\times n}(S)
$$
where each $d_i$ is a divisor of $d_{i+1}$.
The matrix ${\rm diag}(d_1,d_2,\dots)$ is called a \emph{Smith normal form} of the matrix $A$.

If every matrix with coefficients in $S$ admits diagonal reduction, we call $S$ an \emph{elementary divisor ring}.
\end{definition}

\begin{theorem}\cite[Theorem 15.9, p.~188]{Brown}
Any principal ideal ring is an elementary divisor ring.
Consequentially, any finite chain ring, being a principal ideal ring, is an elementary divisor ring.
\end{theorem}

\begin{definition}\label{definition:Smith-normal-form}
Let $k,n \in \mathbb{N}$ and $(n_{1}, n_{2},\dots,n_{r})\in \mathbb{N}^{r}$ be such that $n_{1}+n_{2}+\dots + n_{r}=k\leq n$. Denote by $D_{(n_{1},n_{2},\dots,n_{r})}$ the diagonal matrix in ${\rm M}_{n\times 2n}(R)$ given by
\begin{align*}
D_{(n_{1},n_{2},\dots,n_{r})} &= \operatorname{diag} \{ \underbrace{\pi^{r-1},\dots,\pi^{r-1}}_\text{$n_{1}$ times}, \underbrace{\pi^{r-2},\dots,\pi^{r-2}}_\text{$n_{2}$ times}, \dots, \underbrace{1,\dots,1}_\text{$n_{r}$ times} \} \\
    &=
    \begin{pNiceMatrix}[first-row, first-col]
    &n_1 & n_2 & \dots & n_r & 2n-k & \\
    n_1 & \pi^{r-1}I & 0 & \cdots & 0 & 0\\
    n_2 & 0 & \pi^{r-2}I & \cdots & 0 & 0\\
    \vdots & \vdots & \vdots & \ddots & \vdots & \vdots\\
    n_r & 0 & 0 & \cdots & I & 0\\
    n-k & 0 & 0 & \cdots & 0 & 0
    \end{pNiceMatrix}.
\end{align*}
\end{definition}

\begin{corollary}
\label{corollary:SNF-fcr-normalized}
For $A \in {\rm M}_{n\times 2n}(R)$, there exists a unique tuple $(n_{1},\dots,n_{r})\in \mathbb{N}^{r}$ satisfying $0\leq n_{1}+\dots +n_{r}\leq n$ so that $A$ can be written in the form $UD_{(n_{1},\dots, n_{r})}V$ where $(U,V)\in {\rm GL}_{n}(R)\times {\rm GL}_{2n}(R)$. 
\end{corollary}

\begin{proof}
Because a finite chain ring is an elementary divisor ring, the matrix $A$ has a Smith normal form. Let $D'=\operatorname{diag}\{d_{1},d_{2},\dots,d_{s}\}$ be a Smith normal form of $A$. Since the Smith normal form is unique up to multiplication by units (see \cite[Theorem~15.24, p.~194]{Brown}),
$$
\operatorname{diag}\{\pi^{{\rm val}(d_{1})},\pi^{{\rm val}(d_{2})},\dots,\pi^{{\rm val}(d_{s})}\},
$$
where $0\leq {\rm val}(d_{1}) \leq \cdots \leq {\rm val}(d_{s}) \leq r-1$,
is a Smith normal form of $A$. 
The diagonal matrix 
$$\operatorname{diag}\{\pi^{{\rm val}(d_{1})},\pi^{{\rm val}(d_{2})},\dots,\pi^{{\rm val}(d_{s})}\}$$ 
depends only on $A$ but not $D'$.

Put ${\rm val}(d_{i})=a_{i}$ for $1\leq i\leq s$, and set
\begin{align*}
K' &= \operatorname{diag}\{\pi^{a_{1}},\pi^{a_{2}},\dots,\pi^{a_{s}}\} \\
K &= \operatorname{diag}\{\pi^{a_{s}},\pi^{a_{s-1}},\dots,\pi^{a_{1}}\} .
\end{align*}
Write
\[
I'_{s}=
\begin{pNiceMatrix}[first-row, first-col]
&& & s &&\\
&0 & 0 & \dots & 0 & 1 \\
&0 & 0 & \dots & 1 & 0 \\
s&\vdots & \vdots & \ddots & \vdots & \vdots \\
&0 & 1 & \dots & 0 & 0 \\
&1 & 0 & \dots & 0 & 0
\end{pNiceMatrix}
\]
and for $s\leq t$ 
\[
I'_{s,t}=
\begin{pNiceMatrix}
I'_s & 0 \\
0 & I_{t-s}
\end{pNiceMatrix} .
\]

Since $K'$ is a Smith normal form of $A$, we have $0\leq a_{1}\leq a_{2}\leq \dots \leq a_{s}\leq r-1$ and $A$ can be written in the form $U'K'V'$ where $(U',V')\in {\rm GL}_{n}(R)\times {\rm GL}_{2n}(R)$. 
On the other hand, $A=U'K'V'=U'I'_{s,n}KI'_{s,2n}V'$. Put $U=U'I'_{s,n},V=I'_{s,2n}V'$; plainly $U$ and $V$ are invertible, and $A=UKV$. 
Since $0\leq a_{1}\leq a_{2}\leq \dots \leq a_{s}\leq r-1$, there exists a unique tuple $(n_{1},\dots,n_{r})\in \mathbb{N}^{r}$ such that $0\leq n_{1}+\dots +n_{r}\leq n$ and that $K=D_{(n_{1},\dots,n_{r})}$. 
The corollary follows.
\end{proof}

\section{Counting over finite chain rings}\label{section:counting-fcr}

Throughout this section, we let $R$ denote a finite chain ring, $\pi$ its uniformizer, $R^\times$ its unit group, $k=R/(\pi)$ its residue field, and $q=|k|$ the size of the residue field.

\subsection{Invertible matrices}

In this paragraph we shall count invertible matrices over finite chain rings.
First, we recall a well known fact about the order of a general linear group over a finite field.

\begin{lemma}\label{lemma:invertible-matrices-finite-field}
The cardinality of ${\rm GL}_{n}(\mathbb{F}_q)$ is
$$|{\rm GL}_{n}(\mathbb{F}_q)| = (q^{n}-1)(q^{n}-q)\dots (q^{n}-q^{n-1}).$$
\end{lemma}

\begin{lemma}\label{lemma:invertible-matrices-fcr0}
Suppose that $A, B \in {\rm M}_n(R)$ satisfy $A \equiv B \,({\rm mod} \, (\pi))$, i.e.~$A \in B + \pi \, {\rm M}_n(R)$. 
Then, $A$ is invertible if and only if $B$ is invertible.
\end{lemma}

\begin{proof}
This is immediate on observing that a matrix $T \in {\rm M}_n(R)$ is invertible if and only if $\det(T)\in R^\times$.
\end{proof}

\begin{lemma}\label{lemma:invertible-matrices-fcr}
The cardinality of ${\rm GL}_{n}(R)$ is
\begin{equation}\label{equation:invertible-matrices-fcr-count}
    Q(q,r,n) = q^{(r-1)n^{2}} (q^{n}-1)(q^{n}-q) \dots (q^{n}-q^{n-1}).    
\end{equation}
Consequentially, if ${\rm Z}_{n}(R)$ denotes the set of non-invertible $n \times n$ matrices over $R$, then, as $q\to\infty$, 
\begin{equation}\label{equation:non-invertible-matrices-fcr-count}
|{\rm Z}_{n}(R)| \sim q^{rn^{2}-1} .
\end{equation}
\end{lemma}

\begin{proof}
A matrix $A\in {\rm M}_n(R)$ is invertible if and only if $\det(A)\in R^\times$.
It is easy to see that the reduction homomorphism
\begin{equation*}
{\rm GL}_{n}(R) \to {\rm GL}_{n}(k), 
A \mapsto A \,{\rm mod}\, (\pi)
\end{equation*}
is surjective, and that the size of its kernel is $q^{(r-1)n^{2}}$. 
On the other hand, since the residue field is a finite field $k=\mathbb{F}_q$, say, by Lemma \ref{lemma:invertible-matrices-finite-field} we have 
$$|{\rm GL}_{n}(k)| = (q^{n}-1)(q^{n}-q)\dots (q^{n}-q^{n-1}).$$
Therefore
$$|{\rm GL}_{n}(R)| = q^{(r-1)n^{2}} (q^{n}-1)(q^{n}-q) \dots (q^{n}-q^{n-1}).$$

The second statement is immediate on noting that $|{\rm Z}_{n}(R)| = |{\rm M}_{n}(R)| - |{\rm GL}_{n}(R)|$. The proof is complete.
\end{proof}

For convenience, we make the convention that $Q(q,r,0) = 1$.

\subsection{Matrices with a prescribed Smith normal form}

In this paragraph, we shall prove several counting results concerning matrices over finite chain rings with a prescribed Smith normal form.

\begin{lemma}\label{lemma:SNF-isotropy-group-preparation}
Let $k,n \in \mathbb{N}$ and $(n_{1},\dots,n_{r})\in \mathbb{N}^{r}$. 
Suppose that $n_{1}+\dots + n_{r}=k\leq n$.
Put $D = D_{(n_{1},\dots,n_{r})}$ and
$$ S_{1} = S_1(D) := \{(U,V) \in {\rm GL}_{n}(R)\times {\rm GL}_{2n}(R): UD=DV\}.$$ 
We have
$$|S_{1}| =
    q^{\sum_{1\leq i,j\leq r}(r-\min(i,j))n_{i}n_{j}} \cdot 
    q^{(3n-2k)\sum_{1\leq i\leq r}(r-i)n_{i}} \cdot 
    q^{r(5n^{2}-3nk+k^{2})} \cdot W$$
where
\begin{equation}\label{equation:SNF-isotropy-group-preparation-W}
W = W(q,n,n_{1},\dots,n_{r}) := 
    \frac{Q(q,r,n-k)}{q^{r(n-k)^{2}}} \cdot
    \frac{Q(q,r,2n-k)}{q^{r(2n-k)^{2}}} \cdot
    \prod_{i=1}^{r}\frac{Q(q,r,n_{i})}{q^{rn_{i}^{2}}} .    
\end{equation}

\end{lemma}

\begin{proof}
If $k=0$, then $n_{1}=...=n_{r}=0$ and $D$ is the zero matrix.
Then, by Lemma \ref{lemma:invertible-matrices-fcr},
\begin{equation*}\label{eqn:special-case}
|S_{1}|=|{\rm GL}_{n}(R)|| {\rm GL}_{2n}(R)|=Q(q,r,n)Q(q,r,2n).
\end{equation*}
The lemma is proved for the case $k=0$.

Suppose now that $k \geq 1$. Let $1\leq \alpha_{1}< \alpha_{2}< \dots < \alpha_{t} \leq r$ be the ascending sequence for which $\{n_{\alpha_{1}},n_{\alpha_{2}},\dots,n_{\alpha_{t}}\}$ is the subset of all positive integers of $\{n_{1},\dots ,n_{r}\}$. In particular, $n_{\alpha_{1}}+n_{\alpha_{2}}+\dots+n_{\alpha_{t}}=k$.
By definition,
$$
D=\operatorname{diag} \{\underbrace{\pi^{r-\alpha_{1}},\dots,\pi^{r-\alpha_{1}}}_\text{$n_{\alpha_{1}}$ times}, \underbrace{\pi^{r-\alpha_{2}},\dots,\pi^{r-\alpha_{2}}}_\text{$n_{\alpha_{2}}$ times}, \dots, \underbrace{\pi^{r-\alpha_{t}},\dots,\pi^{r-\alpha_{t}}}_\text{$n_{\alpha_{t}}$ times}\}.
$$

Consider
$$ S_{1}'=S_{1}'(D)= \{{\rm M}_n(R)\times M_{2n}(R):UD=DV\}.$$
If $(U,V)\in S'_1$, we decompose $U,V$ into blocks as follows:
$$
U=
\begin{pNiceMatrix}[first-row, first-col]
&n_{\alpha_{1}} & n_{\alpha_{2}} & \dots & n_{\alpha_{t}} & n-k & \\
n_{\alpha_{1}} & N_{11}^{U} & N_{12}^{U} & \cdots & N_{1t}^{U} & G_{1}^{U}\\
n_{\alpha_{2}} & N_{21}^{U} & N_{22}^{U} & \cdots & N_{2t}^{U} & G_{2}^{U}\\
\vdots & \vdots & \vdots & \ddots & \vdots & \vdots\\
n_{\alpha_{t}} & N_{t1}^{U} & N_{t2}^{U} & \cdots & N_{tt}^{U} & G_{t}^{U}\\
n-k & H_{1}^{U} & H_{2}^{U} & \cdots & H_{t}^{U} & L^{U}
\end{pNiceMatrix},
$$

$$
V=
\begin{pNiceMatrix}[first-row, first-col]
&n_{\alpha_{1}} & n_{\alpha_{2}} & \dots & n_{\alpha_{t}} & 2n-k & \\
n_{\alpha_{1}} & N_{11}^{V} & N_{12}^{V} & \cdots & N_{1t}^{V} & G_{1}^{V}\\
n_{\alpha_{2}} & N_{21}^{V} & N_{22}^{V} & \cdots & N_{2t}^{V} & G_{2}^{V}\\
\vdots & \vdots & \vdots & \ddots & \vdots & \vdots\\
n_{\alpha_{t}} & N_{t1}^{V} & N_{t2}^{V} & \cdots & N_{tt}^{V} & G_{t}^{V}\\
2n-k & H_{1}^{V} & H_{2}^{V} & \cdots & H_{t}^{V} & L^{V}
\end{pNiceMatrix}.
$$

The matrix identity $UD=DV$ unravels as
\begin{align*}
& \begin{pNiceMatrix}
\pi^{r-\alpha_{1}}N_{11}^{U} & \pi^{r-\alpha_{2}}N_{12}^{U} & \cdots & \pi^{r-\alpha_{t}}N_{1t}^{U} & 0\\
\pi^{r-\alpha_{1}}N_{21}^{U} & \pi^{r-\alpha_{2}}N_{22}^{U} & \cdots & \pi^{r-\alpha_{t}}N_{2t}^{U} & 0\\
\vdots & \vdots & \ddots & \vdots & \vdots\\
\pi^{r-\alpha_{1}}N_{t1}^{U} & \pi^{r-\alpha_{2}}N_{t2}^{U} & \cdots & \pi^{r-\alpha_{t}}N_{tt}^{U} & 0\\
\pi^{r-\alpha_{1}}H_{1}^{U} & \pi^{r-\alpha_{2}}H_{2}^{U} & \cdots & \pi^{r-\alpha_{t}}H_{t}^{U} & 0\\
\end{pNiceMatrix} \\
=&
\begin{pNiceMatrix}
\pi^{r-\alpha_{1}}N_{11}^{V} & \pi^{r-\alpha_{1}}N_{12}^{V} & \cdots & \pi^{r-\alpha_{1}}N_{1t}^{V} & \pi^{r-\alpha_{1}}G_{1}^{V}\\
\pi^{r-\alpha_{2}}N_{21}^{V} & \pi^{r-\alpha_{2}}N_{22}^{V} & \cdots & \pi^{r-\alpha_{2}}N_{2t}^{V} & \pi^{r-\alpha_{2}}G_{2}^{V}\\
\vdots & \vdots & \ddots & \vdots & \vdots\\
\pi^{r-\alpha_{t}}N_{t1}^{V} & \pi^{r-\alpha_{t}}N_{t2}^{V} & \cdots & \pi^{r-\alpha_{t}}N_{tt}^{V} & \pi^{r-\alpha_{t}}G_{t}^{V}\\
0&0&\cdots&0&0\\
\end{pNiceMatrix}.
\end{align*}
It follows that
$$
\begin{cases}
\begin{aligned}
\pi^{r-\alpha_{j}}N_{ij}^{U} & =\pi^{r-\alpha_{i}}N_{ij}^{V} \quad &&(1\leq i,j\leq t) \\
\pi^{r-\alpha_{i}}H_{i}^{U} &=0 \quad && (1\leq i\leq t)\\
\pi^{r-\alpha_{i}}G_{i}^{V} & =0 \quad &&(1\leq i\leq t) .
\end{aligned}
\end{cases}
$$
Equivalently, we have
\begin{equation}
\label{proof:SNF-system-of-congruences}
\begin{cases}
\begin{aligned}
N_{ij}^{V} &\in \pi^{\alpha_{i}-\alpha_{j}}N_{ij}^{U} +(\pi^{\alpha_{i}}) \quad 
    && (1\leq j<i\leq t)\\
N_{ij}^{U} &\in \pi^{\alpha_{j}-\alpha_{i}}N_{ij}^{V} +(\pi^{\alpha_{j}}) \quad 
    && (1\leq i<j\leq t)\\
N_{ii}^{U}&\in N_{ii}^{V} +(\pi^{\alpha_{i}}) \quad 
    && (1\leq i\leq t) \\
H_{i}^{U} &\in (\pi^{\alpha_{i}}) \quad 
    && (1\leq i\leq t)\\
G_{i}^{V} &\in (\pi^{\alpha_{i}}) \quad 
    && (1\leq i\leq t).
\end{aligned}
\end{cases}
\end{equation}

From the system of equations \eqref{proof:SNF-system-of-congruences} several observations concerning the invertibility of $U$ and $V$ arise.
We first claim that $U$ is invertible if and only if the submatrices $N_{11}^{U},N_{22}^{U},\dots,N_{tt}^{U}, L^{U}$ are invertible. Indeed, it follows from \eqref{proof:SNF-system-of-congruences} that $N_{ij}^{U} \in (\pi)$ for all $1\leq i<j\leq t$, and that $H_{i}^{U}\in (\pi)$ for all $1\leq i\leq t$. By Lemma \ref{lemma:invertible-matrices-fcr0}, the matrix $U$ is invertible if and only if the matrix
$$
U'=
\begin{pNiceMatrix}
N_{11}^{U} & 0 & \cdots & 0 & G_{1}^{U}\\
N_{21}^{U} & N_{22}^{U} & \cdots & 0 & G_{2}^{U}\\
\vdots & \vdots & \ddots & \vdots & \vdots\\
N_{t1}^{U} & N_{t2}^{U} & \cdots & N_{tt}^{U} & G_{t}^{U}\\
0 & 0 & \cdots & 0 & L^{U}
\end{pNiceMatrix}
$$
is invertible. Since $\det U'=\det N_{11}^{U} \det N_{22}^{U} \dots \det N_{tt}^{U} \det L^{U}$, we have $U'$ is invertible if and only if the submatrices $N_{11}^{U},N_{22}^{U},\dots,N_{tt}^{U}, L^{U}$ are invertible. The first claim is proved.
Our second claim is that $V$ is invertible if and only if the submatrices $N_{11}^{V},N_{22}^{V},\dots,N_{tt}^{V}, L^{V}$ are invertible. This can be shown similarly on noting \eqref{proof:SNF-system-of-congruences} implies that $N_{ij}^{V} \in (\pi)$ for all $1\leq j<i \leq t$ and that $G_{i}^{V} \in (\pi)$ for all $1\leq i\leq t$. 
Our third claim is that $N_{11}^{U},N_{22}^{U},\dots,N_{tt}^{U}$ are invertible if and only if $N_{11}^{V},N_{22}^{V},\dots,N_{tt}^{V}$ are invertible. This is because \eqref{proof:SNF-system-of-congruences} implies that $N_{ii}^{U} \in N_{ii}^{V} +(\pi^{\alpha_{i}})$ for all $1\leq i\leq t$, so by Lemma \ref{lemma:invertible-matrices-fcr0} the third claim follows. 
Combining the above three claims, we conclude that for $(U,V)\in S'_1$, we have $U$ and $V$ are invertible if and only if the submatrices $N_{11}^{U},N_{22}^{U},\dots,N_{tt}^{U}, L^{U}, L^{V}$ are invertible. 

By virtue of Lemma \ref{lemma:invertible-matrices-fcr}, the number of tuples $(N_{11}^{U},N_{22}^{U},\dots,N_{tt}^{U}, L^{U}, L^{V})$ consisting only of invertible matrices is
\begin{equation}\label{proof:SNF-tuple-invertible-matrices}
Q(q,r,n-k) \cdot Q(q,r,2n-k) \cdot \prod_{i=1}^{r}Q(q,r,n_{\alpha_{i}}) .
\end{equation}
Note also that for a fixed tuple $(N_{11}^{U},N_{22}^{U},\dots,N_{tt}^{U}, L^{U}, L^{V})$ of invertible matrices, the number of $(U,V) \in S_{1}$ containing this tuple is equal to the number of $(U,V) \in S_{1}'$ containing the same tuple.

We now count the number of $(U,V) \in S_{1}'$ containing a fixed tuple of invertible matrices
$$(N_{11}^{U},N_{22}^{U},\dots,N_{tt}^{U}, L^{U}, L^{V}).$$ 
Our counting proceeds in two phases.
In the first phase, we carry out five counts, according to the five equations of \eqref{proof:SNF-system-of-congruences}.
First, for $1\leq j<i\leq t$, we count the number of $(N_{ij}^{U}, N_{ij}^{V})$ satisfying the first equation of \eqref{proof:SNF-system-of-congruences}. The number of ways to choose $N_{ij}^{U} \in M_{n_{\alpha_{i}} \times n_{\alpha_{j}}}(R)$ arbitrarily is $q^{rn_{\alpha_{i}}n_{\alpha_{j}}}$. For every chosen $N_{ij}^{U}$, the number of $N_{ij}^{V}$ such that $(N_{ij}^{U}, N_{ij}^{V})$ satisfies the first equation of (\ref{proof:SNF-system-of-congruences}) is $q^{(r-\alpha_{i})n_{\alpha_{i}}n_{\alpha_{j}}}$. Hence, for $1\leq i<j\leq t$, the number of $(N_{ij}^{U}, N_{ij}^{V})$ satisfying the first equation of (\ref{proof:SNF-system-of-congruences}) is $q^{(2r-\alpha_{i})n_{\alpha_{i}}n_{\alpha_{j}}}$.
Second, using the same argument, for $1\leq i<j\leq t$, the number of $(N_{ij}^{U},N_{ij}^{V})$ satisfying the second equation of (\ref{proof:SNF-system-of-congruences}) is $q^{(2r-\alpha_{j})n_{\alpha_{i}}n_{\alpha_{j}}}$.
Third, for $1\leq i\leq t$, the number of $N_{ii}^{V}$ satisfying the third equation of \eqref{proof:SNF-system-of-congruences} is $q^{(r-\alpha_{i})n_{\alpha_{i}}^{2}}$.
The fourth equation of \eqref{proof:SNF-system-of-congruences} implies that the number of $H_{i}^{U}$ is $q^{(n-k)(r-\alpha_{i})n_{\alpha_{i}}}$.
The fifth equation of \eqref{proof:SNF-system-of-congruences} 
implies that the number of $G_{i}^{V}$ is $q^{(2n-k)(r-\alpha_{i})n_{\alpha_{i}}}$.
Therefore, in the first phase, the number of ways to choose submatrices containing the fixed tuple and satisfying \eqref{proof:SNF-system-of-congruences} is
\begin{align*}
& q^{\sum_{1\leq j<i\leq t}(2r-\alpha_{i})n_{\alpha_{i}}n_{\alpha_{j}}} \cdot 
q^{\sum_{1\leq i<j\leq t}(2r-\alpha_{j})n_{\alpha_{i}}n_{\alpha_{j}}} \cdot 
q^{\sum_{1\leq i\leq t}((r-\alpha_{i})n_{\alpha_{i}}^{2})} \\ 
    &\times q^{(n-k)\sum_{1\leq i\leq t}(r-\alpha_{i})n_{\alpha_{i}}} \cdot 
    q^{(2n-k)\sum_{1\leq i\leq t}(r-\alpha_{i})n_{\alpha_{i}}} .
\end{align*}
In the second phase, to fully obtain $U$ and $V$, we also need to choose the remaining submatrices $G_{i}^{U}$ and $H_{i}^{V}$ arbitrarily for $1\leq i\leq t$. There are $q^{r(n-k)n_{\alpha_{i}}}$ matrices $G_{i}^{U}$ and $q^{rn_{\alpha_{i}}(2n-k)}$ matrices $H_{i}^{V}$. Therefore, the number of submatrices $G_{i}^{U}$ and $H_{i}^{V}$ for $1\leq i\leq t$ is
$$
q^{r(3n-2k)\sum_{1\leq i\leq t}n_{\alpha_{i}}}.
$$
Combining the two phases, the number of $(U,V)\in S_{1}'$ containing the fixed tuple is
\begin{align*}
& q^{\sum_{1\leq j<i\leq t}(2r-\alpha_{i})n_{\alpha_{i}}n_{\alpha_{j}}} \cdot 
q^{\sum_{1\leq i<j\leq t}(2r-\alpha_{j})n_{\alpha_{i}}n_{\alpha_{j}}} \cdot 
q^{\sum_{1\leq i\leq t}(r-\alpha_{i})n_{\alpha_{i}}^{2}} \\ 
    &\times q^{(n-k)\sum_{1\leq i\leq t}(r-\alpha_{i})n_{\alpha_{i}}} \cdot 
    q^{(2n-k)\sum_{1\leq i\leq t}(r-\alpha_{i})n_{\alpha_{i}}} \cdot 
    q^{r(3n-2k)\sum_{1\leq i\leq t}n_{\alpha_{i}}} .
\end{align*}
This is also the number of $(U,V) \in S_{1}$ containing the fixed tuple. 

It now follows from \eqref{proof:SNF-tuple-invertible-matrices} that
\begin{align*}
|S_{1}|&=
    q^{\sum_{1\leq j<i\leq t}(2r-\alpha_{i})n_{\alpha_{i}}n_{\alpha_{j}}} \cdot 
    q^{\sum_{1\leq i<j\leq t}(2r-\alpha_{j})n_{\alpha_{i}}n_{\alpha_{j}}} \cdot 
    q^{\sum_{1\leq i\leq t}(r-\alpha_{i})n_{\alpha_{i}}^{2}} \\
&\quad \times
    q^{(n-k)\sum_{1\leq i\leq t}(r-\alpha_{i})n_{\alpha_{i}}} \cdot
    q^{(2n-k)\sum_{1\leq i\leq t}(r-\alpha_{i})n_{\alpha_{i}}} \cdot 
    q^{r(3n-2k)\sum_{1\leq i\leq t}n_{\alpha_{i}}} \\
&\quad \times
    Q(q,r,n-k) \cdot Q(q,r,2n-k) \cdot \prod_{i=1}^{r}Q(q,r,n_{\alpha_{i}}) \\
&= 
    q^{\sum_{1\leq j<i\leq t}(2r-\alpha_{i})n_{\alpha_{i}}n_{\alpha_{j}}} \cdot 
    q^{\sum_{1\leq i<j\leq t}(2r-\alpha_{j})n_{\alpha_{i}}n_{\alpha_{j}}} \cdot 
    q^{\sum_{1\leq i\leq t}(r-\alpha_{i})n_{\alpha_{i}}^{2}} \cdot 
    q^{r\sum_{1\leq i\leq t} n_{\alpha_{i}}^{2}} \\
&\quad\times
    q^{(n-k)\sum_{1\leq i\leq t}(r-\alpha_{i})n_{\alpha_{i}}} \cdot 
    q^{(2n-k)\sum_{1\leq i\leq t}(r-\alpha_{i})n_{\alpha_{i}}} \cdot 
    q^{r(3n-2k)\sum_{1\leq i\leq t}n_{\alpha_{i}}} \cdot 
    q^{r((n-k)^{2}+(2n-k)^{2})}\\
&\quad \times
    \frac{Q(q,r,n-k)}{q^{r(n-k)^{2}}} \cdot
    \frac{Q(q,r,2n-k)}{q^{r(2n-k)^{2}}} \cdot
    \prod_{i=1}^{t}\frac{Q(q,r,n_{\alpha_{i}})}{q^{rn_{\alpha_{i}}^{2}}} \\
&= 
    q^{\sum_{1\leq i,j \leq t} (2r-\alpha_{\min(i,j)})n_{\alpha_{i}}n_{\alpha_{j}}}
\cdot 
    q^{(3n-2k)\sum_{1\leq i\leq t} (r-\alpha_{i})n_{\alpha_{i}}} \cdot 
    q^{r(3n-2k)k} \cdot 
    q^{r(5n^{2}-6nk+2k^{2})} \\
&\quad\times
    \frac{Q(q,r,n-k)}{q^{r(n-k)^{2}}} \cdot
    \frac{Q(q,r,2n-k)}{q^{r(2n-k)^{2}}} \cdot
    \prod_{i=1}^{t}\frac{Q(q,r,n_{\alpha_{i}})}{q^{rn_{\alpha_{i}}^{2}}} \\
&= 
    q^{\sum_{1\leq i,j \leq t} (r-\alpha_{\min(i,j)})n_{\alpha_{i}}n_{\alpha_{j}}}\cdot 
    q^{rk^{2}} \cdot 
    q^{(3n-2k)\sum_{1\leq i\leq t} (r-\alpha_{i})n_{\alpha_{i}}} \cdot 
    q^{r(3n-2k)k} \cdot q^{r(5n^{2}-6nk+2k^{2})} \\
&\quad\times
    \frac{Q(q,r,n-k)}{q^{r(n-k)^{2}}} \cdot
    \frac{Q(q,r,2n-k)}{q^{r(2n-k)^{2}}} \cdot
    \prod_{i=1}^{t}\frac{Q(q,r,n_{\alpha_{i}})}{q^{rn_{\alpha_{i}}^{2}}} \\
&= 
    q^{\sum_{1\leq i,j \leq t} (r-\alpha_{\min(i,j)})n_{\alpha_{i}}n_{\alpha_{j}}}
\cdot 
    q^{(3n-2k)\sum_{1\leq i\leq t} (r-\alpha_{i})n_{\alpha_{i}}} \cdot q^{r(5n^{2}-3nk+k^{2})} \\
&\quad\times
    \frac{Q(q,r,n-k)}{q^{r(n-k)^{2}}} \cdot
    \frac{Q(q,r,2n-k)}{q^{r(2n-k)^{2}}} \cdot
    \prod_{i=1}^{t}\frac{Q(q,r,n_{\alpha_{i}})}{q^{rn_{\alpha_{i}}^{2}}} .
\end{align*}
Adding the numbers $i$ where $n_{i}=0$ to the sums and products, we obtain
$$
|S_{1}| = 
    q^{\sum_{1\leq i,j\leq r}(r-\min(i,j))n_{i}n_{j}} \cdot 
    q^{(3n-2k)\sum_{1\leq i\leq r}(r-i)n_{i}} \cdot 
    q^{r(5n^{2}-3nk+k^{2})} \cdot W .
$$
\end{proof}

\begin{corollary}\label{corollary:SNF-isotropy-group}
Let $k,n \in \mathbb{N}$ and $(n_{1},\dots,n_{r})\in \mathbb{N}^{r}$. 
Suppose that $n_{1}+\dots + n_{r}=k\leq n$.
Put $D = D_{(n_{1},\dots,n_{r})}$ and $S_{1} = S_1(D)$. 
When $n$ is fixed and $q\to \infty$, we have
$$|S_{1}| \sim 
    q^{\sum_{1\leq i,j\leq r}(r-\min(i,j))n_{i}n_{j}} \cdot 
    q^{(3n-2k)\sum_{1\leq i\leq r}(r-i)n_{i}} \cdot q^{r(5n^{2}-3nk+k^{2})}.$$
\end{corollary}

\begin{proof}
It suffices to show that in Lemma \ref{lemma:SNF-isotropy-group-preparation}, the term $W$ defined by \eqref{equation:SNF-isotropy-group-preparation-W} satisfies $W\sim 1$. Suppose that $n$ is fixed. Because $n_{1}+\dots +n_{r} = k \leq n$, the term $W$ have at most $n$ factors of the form $\frac{Q(q,r,n_{i})}{q^{rn_{i}^{2}}}$ where $n_{i}>0$. 
On the other hand, by \eqref{equation:invertible-matrices-fcr-count},
$$
\lim_{q\to \infty} \frac{Q(q,r,n)}{q^{rn^2}} = 1,
$$
whence
$\lim_{q\to \infty} W = 1.$
The corollary is proved.
\end{proof}

\begin{lemma}
\label{lemma:matrices-prescribedSNF}
Let $k,n \in \mathbb{N}$ and $(n_{1},\dots,n_{r})\in \mathbb{N}^{r}$. 
Suppose that $n_{1}+\dots + n_{r}=k\leq n$.
Put $D = D_{(n_{1},\dots,n_{r})}$ and 
$$
S_{2} = S_{2}(D) := 
    \{A \in {\rm M}_{n\times 2n}(R): \exists (U,V) \in {\rm GL}_{n}(R)\times {\rm GL}_{2n}(R)
    \text{ with } A = UDV \}.
$$
When $n$ is fixed and $q\to \infty$, we have
$$
|S_{2}| \sim 
    \frac{q^{r(3nk - k^{2})}}
    {q^{\sum_{1\leq i,j\leq r}(r-\min(i,j))n_{i}n_{j}} \cdot 
    q^{(3n-2k)\sum_{1\leq i\leq r}(r-i)n_{i}}} .
$$
\end{lemma}

\begin{proof}
We observe that the group ${\rm GL}_{n}(R)\times {\rm GL}_{2n}(R)$ acts on $S_1=S_1(D)$ transitively.
The stabilizer of $D$ has cardinality equal to $|S_1|$.
Therefore
$$
|S_{2}| = \frac{|{\rm GL}_{n}(R)\times {\rm GL}_{2n}(R)|}{|S_{1}|}.
$$
It then follows from Lemma \ref{lemma:invertible-matrices-fcr} and Corollary \ref{corollary:SNF-isotropy-group} that
\begin{align*}
|S_{2}| &\sim 
    \frac{Q(q,r,n)Q(q,r,2n)}
    {q^{\sum_{1\leq i,j\leq r}(r-\min(i,j))n_{i}n_{j}} \cdot 
    q^{(3n-2k)\sum_{1\leq i\leq r}(r-i)n_{i}} \cdot q^{r(5n^{2}-3nk+k^{2})}}\\
    &\sim 
    \frac{q^{r(3nk - k^{2})}}
    {q^{\sum_{1\leq i,j\leq r}(r-\min(i,j))n_{i}n_{j}} \cdot 
    q^{(3n-2k)\sum_{1\leq i\leq r}(r-i)n_{i}}} .
\end{align*}
\end{proof}

\subsection{Systems of linear equations}

In this paragraph, we count the number of solutions to a system of linear equations over a finite chain ring.

\begin{lemma}\label{lemma:SNF-image-count}
Let $k,n \in \mathbb{N}$ and $(n_{1},\dots,n_{r})\in \mathbb{N}^{r}$. 
Suppose that $k\leq n$ and that $n_{1}+\dots + n_{r}=k$.
Put $D = D_{(n_{1},\dots,n_{r})}$. 
Suppose that $A\in {\rm M}_{n\times 2n}(R)$ can be written in the form $UDV$ where $(U,V) \in {\rm GL}_{n}(R)\times {\rm GL}_{2n}(R)$. 
Then the number of $b \in R^n$ such that the equation $Ax=b$ has a solution $x \in R^{2n}$ is
$$
q^{\sum_{1\leq i\leq r}in_{i}} .
$$
\end{lemma}

\begin{proof}
The equation $Ax=b$ is equivalent to the equation $D(Vx)=U^{-1}b$. 
Put $y=Vx$ and $b'=U^{-1}b$, so that $Dy=b'$. 
The equation $Ax=b$ has a solution $x \in R^{2n}$ if and only if the equation $Dy=b'$ has a solution $y \in R^{2n}$. 
The number of $b \in R^n$ such that the equation $Dy=U^{-1}b$ has a solution $y \in R^{2n}$ is therefore equal to the number of $b' \in R^n$ such that the equation $Dy=b'$ has a solution $y \in R^{2n}$. 

On writing $y=(y_{1},\dots,y_{2n})\in R^{2n}$ and $b'=(b_{1}',\dots,b_{n}') \in R^{n}$, we see that the equation $Dy=b'$ unravels and becomes 
\begin{equation*}
\begin{cases}
\begin{aligned}
\pi^{r-1}y_{i}&=b_{i}' & \qquad (1\leq i\leq n_1)\\
\pi^{r-2}y_{n_{1}+i}&=b_{n_{1}+i}' & \qquad (1\leq i\leq n_2)\\
&\dots\\
y_{n_{1}+...+n_{r-1}+i}&=b_{n_{1}+...+n_{r-1}+i}' & \qquad (1\leq i\leq n_r)\\
0&=b_{j}' & \qquad (k+1\leq j\leq n).
\end{aligned}
\end{cases}
\end{equation*}
This system of equations has a solution $y \in R^{2n}$ if and only if $b' \in R^{n}$ satisfies
\begin{equation*}
\begin{cases}
\begin{aligned}
b_{i}' &\in (\pi^{r-1}) & \qquad (1\leq i\leq n_1) \\
b_{n_{1}+i}'&\in (\pi^{r-2}) & \qquad (1\leq i\leq n_2) \\
&\dots\\
b_{n_{1}+\dots +n_{r-1}+i}'&\in R & \qquad (1\leq i\leq n_r) \\
b_{j}'&=0 & \qquad (k+1\leq j\leq n) .
\end{aligned}
\end{cases}
\end{equation*}
The number of these solutions $(b_{1}',\dots,b_{n}') \in R^{n}$ is $q^{\sum_{1\leq i\leq r}in_{i}}$. The lemma follows.
\end{proof}

\begin{lemma}\label{lemma:SNF-solution-count}
Let $k,n \in \mathbb{N}$ and $(n_{1},\dots,n_{r})\in \mathbb{N}^{r}$. 
Suppose that $k\leq n$ and that $n_{1}+\dots + n_{r}=k$.
Put $D = D_{(n_{1},\dots,n_{r})}$. 
Let $A\in {\rm M}_{n\times 2n}(R)$ be a matrix that can be written in the form $UDV$ where $(U,V) \in {\rm GL}_{n}(R)\times {\rm GL}_{2n}(R)$. 
Let $b \in R^n$ be a vector so that $Ax=b$ has a solution $x\in R^{2n}$. 
Then the number of solutions $x\in R^{2n}$ to the equation $Ax=b$ is
$$ 
q^{2rn-\sum_{1\leq i\leq r}in_{i}}.
$$
\end{lemma}

\begin{proof}
Let $X_{1}$ be a solution to the equation $Ax=b$. Then, $X$ is a solution to the equation $Ax=b$ if and only if $X-X_{1}$ is a solution to the equation $Ax=0$. Hence the number of solutions $x\in R^{2n}$ to the equation $Ax=b$ is equal to the number of solutions $x\in R^{2n}$ to the equation $Ax=0$.

Since $A=UDV, (U,V)\in {\rm GL}_{n}(R)\times {\rm GL}_{2n}(R)$ and $Ax=0\Leftrightarrow UDVx=0 \Leftrightarrow Dy=0$ where $y=Vx$, the number of solutions $x\in R^{2n}$ to the equation $Ax=0$ is equal to the number of solutions $y\in R^{2n}$ to the equation $Dy=0$. Thus, the number of solutions $x\in R^{2n}$ to the equation $Ax=b$ is equal to the number of solutions $y\in R^{2n}$ to the equation $Dy=0$. We count the number of solutions $y\in R^{2n}$ to the equation $Dy=0$, which is the number of $(y_{1},...,y_{2n})\in R^{2n}$ satisfying the equations
\begin{equation*}
\begin{cases}
\begin{aligned}
\pi^{r-1}y_{i}&=0 & \qquad (1\leq i\leq n_1) \\
\pi^{r-2}y_{n_{1}+i}&=0 & \qquad (1\leq i\leq n_2)\\
&\dots\\
y_{n_{1}+\dots +n_{r-1}+i}&=0 & \qquad (1\leq i\leq n_r) \\
y_{j}&=0 & \qquad (k+1\leq j\leq n) .
\end{aligned}
\end{cases}
\end{equation*}
Thus the number of $y \in R^{2n}$ satisfying the above system of equations is 
$$
q^{\sum_{1\leq i\leq r}n_{i}(r-i)}\cdot q^{r(2n-k)} = q^{2rn-\sum_{1\leq i\leq r}in_{i}}.
$$
\end{proof}

\section{Counting solutions via graphs}\label{section:matrix-equation}

\subsection{Regular graphs}

In this paragraph we shall construct two regular graphs from the matrix equation $AX=C$ and estimate the degrees of these graphs.

\begin{definition}\label{definition:graph-solvable-unsolvable}
Let $k,n \in \mathbb{N}$ and $(n_{1},\dots,n_{r})\in \mathbb{N}^{r}$ satisfying $n_{1}+\dots + n_{r}=k$ and $k\leq n$.
Define $G_{(n_{1},\dots,n_{r})}$ (resp.~$H_{(n_{1},\dots,n_{r})}$) to be the undirected graph whose vertex set is $({\rm M}_n(R))^{3}$ and whose edges are characterized as follows. Two vertices $(a_{1},e_{1},c_{1})$ and $(a_{2},e_{2},c_{2})$ form an edge if the following two conditions are satisfied:
\begin{enumerate}
\item[{\rm (i)}] The equation
\[
\begin{pNiceMatrix}
a_{1}-a_{2} & e_{1}-e_{2}
\end{pNiceMatrix}
x=c_{1}-c_{2}
\]
has a solution (resp.~ has no solution) $x\in {\rm M}_{2n\times n}(R)$.
\item[{\rm (ii)}] The matrix
$\begin{pNiceMatrix}
a_{1}-a_{2} & e_{1}-e_{2}
\end{pNiceMatrix} \in {\rm M}_{2n\times n}(R)$
can be written as $UD_{(n_{1},\dots,n_{r})}V$ for some $(U,V) \in {\rm GL}_{n}(R)\times {\rm GL}_{2n}(R)$.
\end{enumerate}
\end{definition}

\begin{lemma}
\label{lemma:solvable-graph-degree}
Let $k,n \in \mathbb{N}$ and $(n_{1},\dots,n_{r})\in \mathbb{N}^{r}$ satisfying $n_{1}+\dots + n_{r}=k$ and $k\leq n$. 
The graph $G_{(n_{1},\dots,n_{r})}$ is regular and
$$
\deg G_{(n_{1},\dots,n_{r})} \sim 
    \frac{q^{r(3nk-k^{2})}\cdot q^{n\sum_{1\leq i\leq r} in_{i}}}
    {q^{\sum_{1\leq i,j\leq r}(r-\min(i,j))n_{i}n_{j}} \cdot 
        q^{(3n-2k)\sum_{1\leq i\leq r}(r-i)n_{i}}}.
$$
\end{lemma}

\begin{proof}
Observe that two vertices $(a_{1},e_{1},c_{1})$ and $(a_{2},e_{2},c_{2})$ are adjacent to each other if and only if $(0,0,0)$ is adjacent to $(a_{1}-a_{2},e_{1}-e_{2},c_{1}-c_{2})$. Hence $G_{(n_{1},\dots,n_{r})}$ is regular.

Computing the degree of $(0,0,0)$ is tantamount to counting the number of $(a,e,c)$ adjacent to $(0,0,0)$.
We appeal to Lemma \ref{lemma:matrices-prescribedSNF} to compute the number of $(a,e)$ satisfying the condition ${\rm (ii)}$ of Definition \ref{definition:graph-solvable-unsolvable}. With a chosen $(a,e)$, Lemma \ref{lemma:SNF-image-count} gives us the number of ways to choose $c$ satisfying the condition ${\rm (i)}$ of Definition \ref{definition:graph-solvable-unsolvable}. We deduce that
\begin{equation*}
\deg G_{(n_{1},\dots,n_{r})} \sim 
    \frac{q^{r(3nk - k^{2})}}
    {q^{\sum_{1\leq i,j\leq r}(r-\min(i,j))n_{i}n_{j}} \cdot 
        q^{(3n-2k)\sum_{1\leq i\leq r}(r-i)n_{i}}} \cdot 
        q^{n\sum_{1\leq i\leq r}in_{i}}.
\end{equation*}
The lemma is proved.
\end{proof}

\begin{lemma}
\label{lemma:unsolvable-graph-degree}
Let $k,n \in \mathbb{N}$ and $(n_{1},\dots,n_{r})\in \mathbb{N}^{r}$ satisfying $n_{1}+\dots + n_{r}=k$ and $k\leq n$. 
The graph $H_{(n_{1},\dots,n_{r})}$ is regular and
$$\deg H_{(n_{1},\dots,n_{r})} \sim 
    \frac{q^{r(3nk-k^{2})}\cdot (q^{rn^{2}}-q^{n\sum_{1\leq i\leq r} in_{i}})}
    {q^{\sum_{1\leq i,j\leq r}(r-\min(i,j))n_{i}n_{j}} \cdot 
        q^{(3n-2k)\sum_{1\leq i\leq r}(r-i)n_{i}}}.$$
\end{lemma}

\begin{proof}
Using the same argument as in the proof Lemma \ref{lemma:solvable-graph-degree}, we see that $H_{(n_{1},\dots,n_{r})}$ is regular. 

The degree of $H_{(n_{1},\dots,n_{r})}$ equals the number of vertices $(a,e,c)$ that are adjacent to $(0,0,0)$. 
We appeal to Lemma \ref{lemma:matrices-prescribedSNF} to compute the number of $(a,e)$ satisfying the condition ${\rm (ii)}$ of Definition \ref{definition:graph-solvable-unsolvable}. 
Lemma \ref{lemma:SNF-image-count} gives us the number of ways to choose $c$ satisfying the condition ${\rm (i)}$ of Definition \ref{definition:graph-solvable-unsolvable}. We deduce that
\begin{equation*}
\deg H_{(n_{1},\dots,n_{r})} \sim 
    \frac{q^{r(3nk - k^{2})}}{q^{\sum_{1\leq i,j\leq r}(r-\min(i,j))n_{i}n_{j}} \cdot 
    q^{(3n-2k)\sum_{1\leq i\leq r}(r-i)n_{i}}} \cdot 
    (q^{rn^{2}} - q^{n\sum_{1\leq i\leq r}in_{i}}) .
\end{equation*}
The lemma is proved.
\end{proof}

\subsection{A matrix equation}

Given subsets $A,B,C,D,E,F$ of ${\rm M}_n(R)$, let $N(A,B,C,D,E,F)$ denote the number of solutions to the equation
\begin{equation}\label{equation:matrix-solution-count}
    ab+ef=c+d, \quad (a,b,c,d,e,f) \in A\times B\times C\times D\times E\times F.
\end{equation}
The goal of this paragraph is to establish a good estimate for $N(A,B,C,D,E,F)$.

We start with a preparation lemma.

\begin{lemma}\label{lemma:preparation-estimates}
Let
\begin{align*}
P_{1} &= \sum_{k=1}^{n-1} \left( \sum_{i=0}^{r-1} \frac{q^{r(2n^{2}+3nk-k^{2})}\cdot C^{i}_{k+i-1}}{q^{i(1+n)}} \right)\\
P_{2} &= \sum_{i=1}^{r-1} \frac{q^{4rn^2}\cdot C^{i}_{n+i-1}}{q^{i(1+n)}} .
\end{align*}
Then there exist constants $C_1(n)$ and $C_2(n)$ such that
\begin{align*}
    P_{1} &\leq C_1(n) q^{r(4n^{2}-n-1)} \\
    P_{2} &\leq C_2(n) q^{4rn^2-n-1} .
\end{align*}
\end{lemma}

\begin{proof}
For $1\leq k\leq n$ and $i\in \mathbb{N}^{*}$, 
\begin{equation}
C^{i}_{k+i-1} = \frac{(i+1)(i+2)\dots (i+k-1)}{(k-1)!} 
    \leq (i+k-1)^{k-1} 
    \leq (ik)^{k-1}. \label{proof:SNF-system-of-congruences0}
\end{equation}
We have, by \eqref{proof:SNF-system-of-congruences0},
\begin{align*}
P_{1} &= \sum_{k=1}^{n-1} q^{r(2n^{2}+3nk-k^{2})}  
        \left( 1 + \sum_{i=1}^{r-1} \frac{C^{i}_{k+i-1}}{q^{i(1+n)}} \right) \\
    &\leq \sum_{k=1}^{n-1} q^{r(2n^{2}+3nk-k^{2})} 
        \left( 1 + \sum_{i=1}^{r-1} \frac{(ik)^{k-1}}{q^{i(1+n)}} \right) \\
    &\ll_n \sum_{k=1}^{n-1} q^{r(2n^{2}+3nk-k^{2})} 
    \leq nq^{r(2n^{2}+3n(n-1)-(n-1)^{2})} \\
    &\ll q^{r(4n^{2}-n-1)} ,
\end{align*}
whence the upper bound for $P_1$.

We next show the upper bound for $P_2$. 
By \eqref{proof:SNF-system-of-congruences0}, we have
\begin{align*}
P_{2} &= q^{4rn^2-n-1} 
        \sum_{i=1}^{r-1} \frac{C^{i}_{n+i-1}}{q^{(i-1)(1+n)}} 
    \leq q^{4rn^2-n-1} 
        \sum_{i=1}^{r-1} \frac{(in)^{n-1}}{q^{(i-1)(1+n)}} \\
    &\ll_n q^{4rn^2-n-1} 
        \sum_{j=0}^{\infty} \frac{j^{n}}{2^{jn}} 
    \ll_n q^{4rn^{2}-n-1} ,
\end{align*}
whence the upper bound for $P_2$.
\end{proof}

We are in a position to state the key estimate on the number of matrix solutions to the equation
\eqref{equation:matrix-solution-count}.

\begin{proposition}
\label{proposition:count-matrix-solutions}
For $n\in \mathbb{N}$, there exists $C(n)>0$ such that the following holds. If $A,B,C,D,E,F$ are  subsets of ${\rm M}_n(R)$, then
$$
\left| N(A,B,C,D,E,F) - \frac{|A||B||C||D||E||F|}{q^{rn^{2}}} \right| 
    \leq C(n) q^{2rn^2-\frac{n+1}{2}}\sqrt{|A||B||C||D||E||F|} .
$$
\end{proposition}

\begin{proof}
Consider the undirected bipartite graph $Z=(X \cup Y,E)$ where $X=Y=({\rm M}_n(R))^{3}$ and the vertex set is defined as follows. Two vertices $(a,e,c)\in X$ and $(b,f,d)\in Y$ form an edge if and only if $ab+ef=c+d$. Plainly
$$
|X|=|Y|=|({\rm M}_n(R))^{3}|=q^{3rn^{2}}.
$$
It is apparent that $\deg (X) = \deg (Y) = q^{2rn^2}$, and so
$$
\frac{\deg (X)}{|Y|}=\frac{1}{q^{rn^{2}}}.
$$
In fact, this is because for each vertex $(a,e,c)\in X$ and each pair $(b,f)\in ({\rm M}_n(R))^{2}$, the matrix $d=ab+ef-c$ is uniquely determined. 
The adjacency matrix $A_{Z}$ has the form
$$
A_{Z}=
\begin{pNiceMatrix}
0 & N \\
N^{T} & 0
\end{pNiceMatrix}.
$$
The $|X|\times |Y|$ matrix $N$ is characterized by the following property: $N_{xy}=1$ if and only if there is an edge between $x \in X$ and $y\in Y$ and $N_{xy}=0$ otherwise. 

For two vertices $(a_{1},e_{1},c_{1})$ and $(a_{2},e_{2},c_{2})$ in $X$, let us count the number of their common neighbors, i.e.~the number of solutions $(b,f,d)\in Y$ to the system of equations
\begin{equation}\label{proof:count-common-neighbors1}
\begin{cases}
a_{1}b+e_{1}f=c_{1}+d \\
a_{2}b+e_{2}f=c_{2}+d .
\end{cases}
\end{equation}
It follows that
\begin{equation}\label{proof:count-common-neighbors2}
\begin{pNiceMatrix}
a_{1}-a_{2} & e_{1}-e_{2}
\end{pNiceMatrix}
\begin{pNiceMatrix}
b\\
f
\end{pNiceMatrix}
=
c_{1}-c_{2} .
\end{equation}
The set of solutions to \eqref{proof:count-common-neighbors2} is in bijection to the set of solutions to \eqref{proof:count-common-neighbors1}
via the map sending $\begin{pNiceMatrix} b\\ f \end{pNiceMatrix}$ to
$(b,f,a_{1}b+e_{1}f-c)$. 

By virtue of Corollary \ref{corollary:SNF-fcr-normalized}, there exists a unique tuple $(n_{1},\dots,n_{r})\in \mathbb{N}^{r}$ such that $0\leq n_{1}+\dots +n_{r}=k\leq n$, and that the matrix $\begin{pNiceMatrix} a_{1}-a_{2}& e_{1}-e_{2} \end{pNiceMatrix}$ can be written in the form $UD_{(n_{1},\dots, n_{r})}V$ where $(U,V)\in {\rm GL}_{n}(R)\times {\rm GL}_{2n}(R)$. It now follows from Lemma \ref{lemma:SNF-solution-count} that, if equation \eqref{proof:count-common-neighbors2}
has a solution, then the number of solutions is
$$
    (q^{2rn-\sum_{1\leq i\leq r}in_{i}})^{n}=q^{2rn^{2}-n\sum_{1\leq i\leq r}in_{i}}.
$$

Recall that the graphs $G_{(n_{1},\dots,n_{r})}$ and $H_{(n_{1},\dots,n_{r})}$ are given by Definition \ref{definition:graph-solvable-unsolvable}.
Let $E_{(n_{1},\dots,n_{r})}$ and $F_{(n_{1},\dots,n_{r})}$ be the adjacency matrices of $G_{(n_{1},\dots,n_{r})}$ and $H_{(n_{1},\dots,n_{r})}$ respectively.
We have:

\begin{align}
NN^{T}&=\deg (X) I + 
    \sum_{k=1}^{n} 
        \sum_{\substack{(n_{1},\dots,n_{r})\in \mathbb{N}^{r}\\ n_{1}+\dots +n_{r}=k}} q^{2rn^{2}-n\sum_{1\leq i\leq r}in_{i}} \cdot E_{(n_{1},\dots,n_{r})} \nonumber \\
    &= q^{rn^{2}}J + (\deg (X)-q^{rn^{2}})I 
        \nonumber \\
    & + \sum_{k=1}^{n} 
        \sum_{\substack{(n_{1},\dots,n_{r})\in \mathbb{N}^{r}\\ n_{1}+\dots +n_{r}=k}}  \left( q^{2rn^{2}-n\sum_{1\leq i\leq r}in_{i}} - q^{rn^{2}}\right) \cdot E_{(n_{1},\dots,n_{r})} 
        \nonumber \\
    & - \sum_{k=0}^{n} \sum_{\substack{(n_{1},\dots,n_{r})\in \mathbb{N}^{r}\\ n_{1}+\dots +n_{r}=k}} q^{rn^{2}} \cdot F_{(n_{1},\dots,n_{r})} 
        \label{proof:eigenvalue-square}.
\end{align}

Let $w^{(3)}=(u_{1},\dots,u_{|X|},v_{1},\dots,v_{|Y|})$ be an eigenvector of $A_{G}$ corresponding to the eigenvalue $\lambda_{3}$. In view of Lemma \ref{lemma:prelim-biregular-eval}, the vector $(u_{1},\dots,u_{|X|})^{T}$ is an eigenvector of $NN^{T}$ corresponding to the eigenvalue $\lambda_{3}^{2}$. It follows from \eqref{proof:eigenvalue-square} that
\begin{equation}
\label{proof:eigenvalue-square1}
(\lambda_{3}^{2}-\deg(X) +q^{rn^{2}})(u_{1},\dots,u_{|X|})^{T}=(A_{1}-A_{2})(u_{1},\dots,u_{|X|})^{T}
\end{equation}
where
\begin{align*}
    A_{1} &:= \sum_{k=1}^{n}
        \sum_{\substack{(n_{1},\dots,n_{r})\in \mathbb{N}^{r}\\ n_{1}+\dots +n_{r}=k}} 
        \left( q^{2rn^{2}-n\sum_{1\leq i\leq r}in_{i}} - q^{rn^{2}}\right) \cdot E_{(n_{1},\dots,n_{r})}\\
    A_{2} &:= \sum_{k=0}^{n} 
        \sum_{\substack{(n_{1},\dots,n_{r})\in \mathbb{N}^{r}\\ n_{1}+\dots +n_{r}=k}} 
        q^{rn^{2}} \cdot F_{(n_{1},\dots,n_{r})} .
\end{align*}

Let $\lambda_{E}$ be an eigenvalue of $E_{(n_{1},\dots,n_{r})}$. Since $G_{(n_{1},\dots,n_{r})}$ is regular, Lemma \ref{lemma:prelim-regular-graph-largest-eval} implies that $|\lambda_{E}|\leq \deg(G_{(n_{1},\dots,n_{r})})$. We then apply Lemma \ref{lemma:solvable-graph-degree} to deduce that
$$
|\lambda_{E}| \ll 
    \frac{q^{r(3nk-k^{2})}\cdot q^{n\sum_{1\leq i\leq r} in_{i}}}
    {q^{\sum_{1\leq i,j\leq r}(r-\min(i,j))n_{i}n_{j}} \cdot 
        q^{(3n-2k)\sum_{1\leq i\leq r}(r-i)n_{i}}} 
    =: M^{E}_{(n_{1},\dots,n_{r})}.
$$
Similarly, if $\lambda_{F}$ is an eigenvalue of $F_{(n_{1},\dots,n_{r})}$, then
$$
|\lambda_{F}| \ll 
    \frac{q^{r(3nk-k^{2})}\cdot (q^{rn^{2}}-q^{n\sum_{1\leq i\leq r} in_{i}})}
    {q^{\sum_{1\leq i,j\leq r}(r-\min(i,j))n_{i}n_{j}} \cdot 
        q^{(3n-2k)\sum_{1\leq i\leq r}(r-i)n_{i}}} 
    =: M^{V}_{(n_{1},\dots,n_{r})}. 
$$
Now, if $\lambda$ is an eigenvalue of $A_{1}-A_{2}$, then by (\ref{proof:eigenvalue-square1}) and Lemma \ref{lemma:prelim-sum-symmetric-matrices-eigenvalue} we have
\begin{align}
|\lambda| \ll 
    &\sum_{k=1}^{n} 
    \sum_{\substack{(n_{1},\dots,n_{r})\in \mathbb{N}^{r}\\ n_{1}+\dots +n_{r}=k}} 
    \left( q^{2rn^{2}-n\sum_{1\leq i\leq r}in_{i}} - q^{rn^{2}}\right)  M^{E}_{(n_{1},\dots,n_{r})} \nonumber \\
    &+ \sum_{k=0}^{n} 
    \sum_{\substack{(n_{1},\dots,n_{r})\in \mathbb{N}^{r}\\ n_{1}+\dots +n_{r}=k}}  q^{rn^{2}}  M^{V}_{(n_{1},\dots,n_{r})}
    \label{proof:eigenvalue-square2} .
\end{align}
We then consider two cases. First, for $k = n$ and $(n_{1},\dots,n_{r-1},n_{r})=(0,\dots,0,n)$ we have
$$
\begin{cases}
\left( q^{2rn^{2}-n\sum_{1\leq i\leq r}in_{i}} - q^{rn^{2}}\right) 
M^{E}_{(n_{1},\dots,n_{r})} = 0 \cdot M^{E}_{(n_{1},\dots,n_{r})} = 0\\
q^{rn^{2}} M^{V}_{(n_{1},\dots,n_{r})}=q^{rn^{2}}\cdot 0 = 0 .
\end{cases}
$$
Secondly, in the situation $k \neq n$ or the situation $k = n$ and $(n_{1},\dots,n_{r-1},n_{r}) \neq (0,\dots,0,n)$, 
\begin{align*}
    \left( q^{2rn^{2}-n\sum_{1\leq i\leq r}in_{i}} - q^{rn^{2}}\right) 
         M^{E}_{(n_{1},\dots,n_{r})}
    &\qquad \leq q^{2rn^{2}-n\sum_{1\leq i\leq r}in_{i}} 
        M^{E}_{(n_{1},\dots,n_{r})}\\
    &\qquad = \frac{q^{r(2n^{2}+3nk-k^{2})}}
        {q^{\sum_{1\leq i,j\leq r}(r-\min(i,j))n_{i}n_{j}}  
        q^{(3n-2k)\sum_{1\leq i\leq r}(r-i)n_{i}}}
\end{align*}
and
\begin{align*}
    q^{rn^{2}} M^{V}_{(n_{1},\dots,n_{r})} 
    \leq \frac{q^{r(2n^{2}+3nk-k^{2})}}
        {q^{\sum_{1\leq i,j\leq r}(r-\min(i,j))n_{i}n_{j}} \cdot 
        q^{(3n-2k)\sum_{1\leq i\leq r}(r-i)n_{i}}} .
\end{align*}
Combining with (\ref{proof:eigenvalue-square2}), we obtain
\begin{equation}
\label{proof:eigenvalue-square3}
|\lambda|\ll q^{2rn^2} + 2 \sum^{n}_{k=1}\left( \sum_{\substack{(n_{1},\dots,n_{r})\in \mathbb{N}^{r}\\ n_{1}+\dots +n_{r}=k\\ n_{r}\neq n}} \frac{q^{r(2n^{2}+3nk-k^{2})}}{T_{(n_{1},\dots,n_{r})}} \right)
\end{equation}
where
$$
T_{(n_{1},\dots,n_{r})} := q^{\sum_{1\leq i,j\leq r}(r-\min(i,j))n_{i}n_{j}} \cdot q^{(3n-2k)\sum_{1\leq i\leq r}(r-i)n_{i}}.
$$

For $1\leq k\leq n, 1\leq i\leq r$, let
\begin{equation*}
M_{i,k} = \{(n_{1},\dots,n_{r}) \in \mathbb{N}^{r}: n_{1}+\dots +n_{r}=k, 
n_{1}=\dots =n_{i-1}=0,  n_{i}\geq 1\}.
\end{equation*}
Plainly $M_{i,k} \cap M_{j,k} = \emptyset$ for all $1\leq i<j\leq r$ and
\begin{equation}
\label{proof:eigenvalue-square4}
\bigcup_{1\leq i\leq r} M_{i,k} = 
    \{ (n_{1},\dots, n_{r})\in \mathbb{N}^{r}: n_{1}+\dots +n_{r}=k \}.
\end{equation}
Also, for $1\leq i\leq r$ and $1\leq k\leq n$,
\begin{equation}
\label{proof:eigenvalue-square5}
|M_{i,k}|=C^{r-i}_{k+r-i-1} .
\end{equation}
Now let $h_{j} = \sum_{i=1}^{j} n_{i}$ where $1\leq j\leq r$. It is easy to see that
$$
T_{(n_{1},\dots,n_{r})} = q^{\sum_{1\leq j \leq r-1} (h_{j}^{2}+(3n-2k)h_{j})}.
$$
If $(n_{1},\dots, n_{r})\in M_{i,k}$ then $h_{j}\geq 1$ for all $i\leq j\leq r$. Hence
\begin{equation*}
T_{(n_{1},\dots,n_{r})}=q^{\sum_{1\leq j \leq r-1} (h_{j}^{2}+(3n-2k)h_{j})}\geq q^{(r-i)(1+3n-2k)} \geq q^{(r-i)(1+n)}.
\end{equation*}
Therefore, for $(n_{1},\dots, n_{r})\in M_{i,k}$,
\begin{equation}\label{proof:eigenvalue-square6}
T_{(n_{1},\dots,n_{r})}\geq q^{(r-i)(1+n)}.
\end{equation}

From \eqref{proof:eigenvalue-square3}, \eqref{proof:eigenvalue-square4}, \eqref{proof:eigenvalue-square5}, \eqref{proof:eigenvalue-square6} we infer that
\begin{align}
|\lambda| 
    &\ll q^{2rn^2} + 2 \sum_{k=1}^{n} \sum_{i=1}^{r}  
        \sum_{\substack{(n_{1},\dots, n_{r}) \in M_{i,k} \\ n_{r} \neq n}} \frac{q^{r(2n^{2}+3nk-k^{2})}}{T_{(n_{1},\dots,n_{r})}}  
        \nonumber  \\
    &= q^{2rn^2} + 2\sum_{k=1}^{n-1} \sum_{i=1}^{r} 
        \sum_{(n_{1},\dots, n_{r}) \in M_{i,k}} \frac{q^{r(2n^{2}+3nk-k^{2})}}{T_{(n_{1},\dots,n_{r})}}  
     +2\sum^{r-1}_{i=1}  \sum_{(n_{1},\dots, n_{r}) \in M_{i,n}} \frac{q^{4rn^2}}{T_{(n_{1},\dots,n_{r})}} 
     \nonumber \\
    & \leq q^{2rn^2} 
        +2\sum_{k=1}^{n-1}  \sum_{i=1}^{r} \frac{q^{r(2n^{2}+3nk-k^{2})}\cdot C^{r-i}_{k+r-i-1}}{q^{(r-i)(1+n)}} 
        +2\sum_{i=1}^{r-1} \frac{q^{4rn^2}\cdot C^{r-i}_{n+r-i-1}}{q^{(r-i)(1+n)}} 
        \nonumber \\
    & = q^{2rn^2} 
        +2\sum_{k=1}^{n-1}  \sum_{i=0}^{r-1} \frac{q^{r(2n^{2}+3nk-k^{2})}\cdot C^{i}_{k+i-1}}{q^{i(1+n)}} 
        +2\sum_{i=1}^{r-1} \frac{q^{4rn^2}\cdot C^{i}_{n+i-1}}{q^{i(1+n)}} 
        \nonumber \\
    &=: q^{2rn^2}+2A_{3}+2A_{4} 
    \label{proof:eigenvalue-square7}.
\end{align}

We now apply Lemma \ref{lemma:preparation-estimates} to derive that
$A_{3}\ll q^{r(4n^{2}-n-1)}$ và $A_{4}\ll q^{4rn^2-n-1}$.
Combining with \eqref{proof:eigenvalue-square7}, we arrive at
$$|\lambda|\ll q^{4rn^2-n-1}$$
where $\lambda$ is an eigenvalue of $A_{1}-A_{2}$. 
By \eqref{proof:eigenvalue-square1} we have
$$
|\lambda_{3}^{2}-\deg (X)+q^{rn^{2}}|\ll q^{4rn^2-n-1}.
$$
Because $\deg (X) = q^{2rn^2}$, this implies that
$$
|\lambda_{3}|\ll q^{2rn^2-\frac{n+1}{2}}.
$$

Finally, observe that for $A,B,C,D,E,F \subset {\rm M}_n(R)$, we have $A\times E\times C \subset X$, $B\times F\times D\subset Y$ and $N(A,B,C,D,E,F) = e(A\times E\times C, B\times D\times F)$. We can now apply Lemma \ref{lemma:expander-mixing} to obtain
\begin{align*}
    \left| N(A,B,C,D,E,F) - \frac{|A||B||C||D||E||F|}{q^{rn^{2}}} \right|
    &\leq |\lambda_{3}|\sqrt{|A\times E\times C||B\times D\times F|}\\
    &\leq C(n) q^{2rn^2-\frac{n+1}{2}}\sqrt{|A||B||C||D||E||F|} .
\end{align*}
The proposition is proved.
\end{proof}

\section{Expander and sum-product on matrix algebras}\label{section:proofs-of main-theorems}

This section contains the proofs of the main results, namely Theorems \ref{theorem-result:expander} and \ref{theorem-result:sumproduct}. 
The proofs of these theorems are applications of Proposition \ref{proposition:count-matrix-solutions}.
We follow the arguments of \cite{XG} and reproduce them here for the sake of completeness.

\begin{proof}[Proof of Theorem \ref{theorem-result:expander}]
For $\lambda \in A+BC$, put
$$
t(\lambda) = |\{ (a,b,c)\in A\times B\times C: a+bc=\lambda \}| .
$$
Using Cauchy-Schwarz inequality, we have
$$
(|A||B||C|)^{2} = \left( \sum_{\lambda \in A+BC} t(\lambda) \right)^{2} \leq |A+BC|\sum_{\lambda \in A+BC} t(\lambda)^{2}.
$$
Note that
$$
\sum_{\lambda \in A+BC} t(\lambda)^{2}=N(B,C,A,-A,-B,C).
$$
On applying Proposition \ref{proposition:count-matrix-solutions} we derive
$$
\frac{(|A||B||C|)^{2}}{|A+BC|}\leq N(B,C,A,-A,-B,C)\ll \frac{|A|^{2}|B|^{2}|C|^{2}}{q^{rn^{2}}}+q^{2rn^2-\frac{n+1}{2}}|A||B||C| .
$$
Thus
$$
|A+BC|\gg \min \left( q^{rn^{2}}, \frac{|A||B||C|}{q^{2rn^2-\frac{n+1}{2}}} \right).
$$
The proof is complete.
\end{proof}

For $A,B \in {\rm M}_n(R)$, define the additive energy
$$
E_{+}(A,B)=|\{ (a_{1},a_{2},b_{1},b_{2})\in A^{2}\times B^{2}:a_{1}+b_{1}=a_{2}+b_{2} \}| .
$$

\begin{lemma}
\label{lemma:additive-energy}
For $A,B\subset {\rm M}_n(R)$ and $C\subset {\rm GL}_{n}(R)$, we have
$$
\frac{(|A||B|)^{2}}{|A+B|}\leq E_{+}(A,B)\ll \frac{|BC|^{2}|A|^{2}}{q^{rn^{2}}}+q^{2rn^2-\frac{n+1}{2}} \frac{|BC||A|}{|C|} .
$$
\end{lemma}

\begin{proof}
We have
\begin{align*}
E_{+}(A,B)&=|\{ (a_{1},a_{2},b_{1},b_{2})\in A^{2}\times B^{2}:a_{1}+b_{1}=a_{2}+b_{2} \}|\\
&= |C|^{-2}|\{ (a_{1},a_{2},b_{1},b_{2},c_{1},c_{2})\in A^{2}\times B^{2}\times C^{2}: a_{1}+b_{1}c_{1}c_{1}^{-1} = a_{2}+b_{2}c_{2}c_{2}^{-1} \}|\\
&\leq |C|^{-2}|\{ (a_{1},a_{2},s_{1},s_{2},t_{1},t_{2})\in A^{2}\times (BC)^{2}\times (C^{-1})^{2}: a_{1}+s_{1}t_{1}= a_{2}+s_{2}t_{2} \}|\\
&= |C|^{-2}| N(BC,C^{-1},A,-A,BC,C^{-1}) .
\end{align*}
On applying Proposition \ref{proposition:count-matrix-solutions}, we see that
\begin{align*}
E_{+}(A,B)&\leq |C|^{-2}| N(BC,C^{-1},A,-A,BC,C^{-1})\\
&\ll |C|^{-2}|\left( \frac{|BC|^{2}|C|^{2}|A|^{2}}{q^{rn^{2}}} + q^{2rn^2-\frac{n+1}{2}} |BC||C||A| \right)\\
&= \frac{|BC|^{2}|A|^{2}}{q^{rn^{2}}} + q^{2rn^2-\frac{n+1}{2}}\frac{|BC||A|}{|C|} .
\end{align*}

We next show that $\frac{(|A||B|)^{2}}{|A+B|}\leq E_{+}(A,B)$. For $\lambda \in A+B$, put
$$
t_{A+B}(\lambda)=|\{(a,b)\in A\times B:a+b=\lambda\}| .
$$
Using Cauchy-Schwarz inequality,
\begin{align*}
(|A||B|)^{2}&=\left(\sum_{\lambda \in A+B} t_{A+B}(\lambda)\right)^{2}
\leq |A+B|\sum_{\lambda \in A+B} t_{A+B}(\lambda)^{2} 
= |A+B|E_{+}(A,B) .
\end{align*}
Thus
$$
\frac{(|A||B|)^{2}}{|A+B|}\leq E_{+}(A,B).
$$
The lemma is proved.
\end{proof}

\begin{proof}[Proof of Theorem \ref{theorem-result:sumproduct}]
Let $Z_{n}(R)={\rm M}_n(R)\backslash {\rm GL}_n(R)$ denote the set of singular matrices in ${\rm M}_n(R)$. Since $|A|\geq C(n)q^{rn^{2}-1}$ and $|Z_{n}(R)|\sim q^{rn^{2}-1}$, we choose $C(n)$ such that $|A|>2|Z_{n}(R)|$. Then, $|A\cap {\rm GL}_{n}(R)|\geq |A|/2$, so we can assume that  $A\subset {\rm GL}_{n}(R)$. On applying Lemma \ref{lemma:additive-energy} with $A=B=C$, we have
$$
\frac{|A|^{4}}{|A+A|}\leq E_{+}(A,A)
    \ll \frac{|AA|^{2}|A|^{2}}{q^{rn^{2}}} + q^{2rn^2-\frac{n+1}{2}}|AA| .$$
Thus
$$
\max( |A+A|,|AA| )
    \gg \min \left( \frac{|A|^{2}}{q^{rn^{2}-\frac{n+1}{4}}} , q^{\frac{rn^{2}}{3}}|A|^{\frac{2}{3}} \right) .
$$
The proof is complete.
\end{proof}

\bibliographystyle{amsplain}

\bibliography{expander}

\end{document}